# An approach for designing a surface pencil through a given asymptotic curve

Fatma GÜLER, Gülnur ŞAFFAK ATALAY, Ergin BAYRAM, Emin KASAP

Ondokuz Mayıs University, Faculty of Arts and Sciences, Mathematics Department
f.guler@omu.edu.tr, g.saffak@omu.edu.tr, erginbayram@yahoo.com, kasape@omu.edu.tr

## ABSTRACT

Surfaces and curves play an important role in geometric design. In recent years, problem of finding a surface passing through a given curve has attracted much interest. In the present paper, we propose a new method to construct a surface interpolating a given curve as the asymptotic curve of it. Also, we analyze the conditions when the resulting surface is a ruled surface. Furthermore, we prove that there exists no developable surface possessing a given curve as an asymptotic curve except plane. Finally, we illustrate this method by presenting some examples.

## 1. Introduction

A rotation minimizing frame (RMF) {T,U,V} of a space curve contains the curve tangent T and the normal plane vectors U, V which show no instantaneous rotation about T. Because of their minimum twist RMFs are very interesting in computer graphics, including free-form deformation with curve constraints [1 - 6], sweep surface modeling [7 - 10], modeling of generalized cylinders and tree branches [11 - 14], visualization of streamlines and tubes [15 - 17], simulation of ropes and strings [18], and motion design and control [19].

There are infinitely many frames on a given space curve [20]. One can produce other frames from an existing one by controlling the orientation of the frame vectors U and V in the normal plane of the curve. In differential geometry the most familiar frame is Frenet frame {T, N, B}, where T is the curve tangent; N is the principal normal vector and $B = T \times N$ is the binormal vector (see [21] for details). Beside of its fame, the Frenet frame is not an RMF and it is unsuitable for specifying the orientation of a rigid body along a given curve in applications such as motion planning, animation, geometric design, and robotics, since it incurs "unnecessary" rotation of the body [22]. Furthermore, Frenet frame is undefined if the curvature vanishes.

An asymptotic curve is one of the most significant surface curves. A surface curve is classified as asymptotic if its velocity is always an asymptotic direction, that is, the direction in which the normal curvature is zero. Furthermore, for a curve on a surface M to be asymptotic, its acceleration must always be tangent to M. The surface is not bending away from its tangent plane in an asymptotic direction.

Asymptotic curves are encountered in astronomy, astrophysics and architectural CAD. In order to find a set of escaping orbits of stars in a stellar system, it is necessary to find asymptotic curves of the Lyapunov orbits because any small outward deviation from an asymptotic orbit will lead a star to escape from the system. Contopoulos [23] considered the asymptotic orbits of mainly unstable orbits, with a particular emphasis on the Lyapunov orbits, and found sets of escaping orbits with initial conditions on asymptotic curves. Efthymiopoulos et al. [24] concluded that the diffusion of any chaotic orbit inside the cantorus follows essentially the same path defined by the unstable asymptotic curves that emanate from unstable periodic orbits inside the cantorus. Flöry and Pottmann [25] addressed challenges in the realization of free-form architecture and complex shapes in general with the technical advantages of ruled surfaces. They proposed a geometry-processing framework to approximate a given shape by one or more strips of ruled surfaces. In that work, they used asymptotic curves obtained by careful investigation and constructed an initial ruled surface by aligning the rulings with asymptotic curves; they also discussed how the shape of this initial approximation can be modified to optimally fit a given target shape.

Most previous work on asymptotic curves has focused on how to find such curves on a given surface. In practice, the more relevant problem is how to construct surfaces that contain a given spatial curve as a common asymptotic curve. We use an example in architectural design to illustrate potential applications of this problem. Figure 1 shows an architectural model that is represented by multiple strips of ruled surfaces. Each point on each ruled surface falls on a straight line laying entirely on that surface. It is obvious that the surface does not curve in the ruling direction. Such a tangential direction is said to be of vanishing normal curvature and is called an asymptotic direction. Excluding developable surfaces, ruled surfaces are known to have a negative Gaussian curvature. Hence, they exhibit another tangential direction of zero normal curvature different from the ruling direction at any point [26]. The asymptotic curves pass through the directions of vanishing normal curvature on a surface. In negatively curved, namely locally saddle-shaped, regions, there are two asymptotic directions for every point on the surface. If the curvature of a set of asymptotic curves in a region is sufficiently small, there is a good possibility that the region can be approximated by a ruled surface.

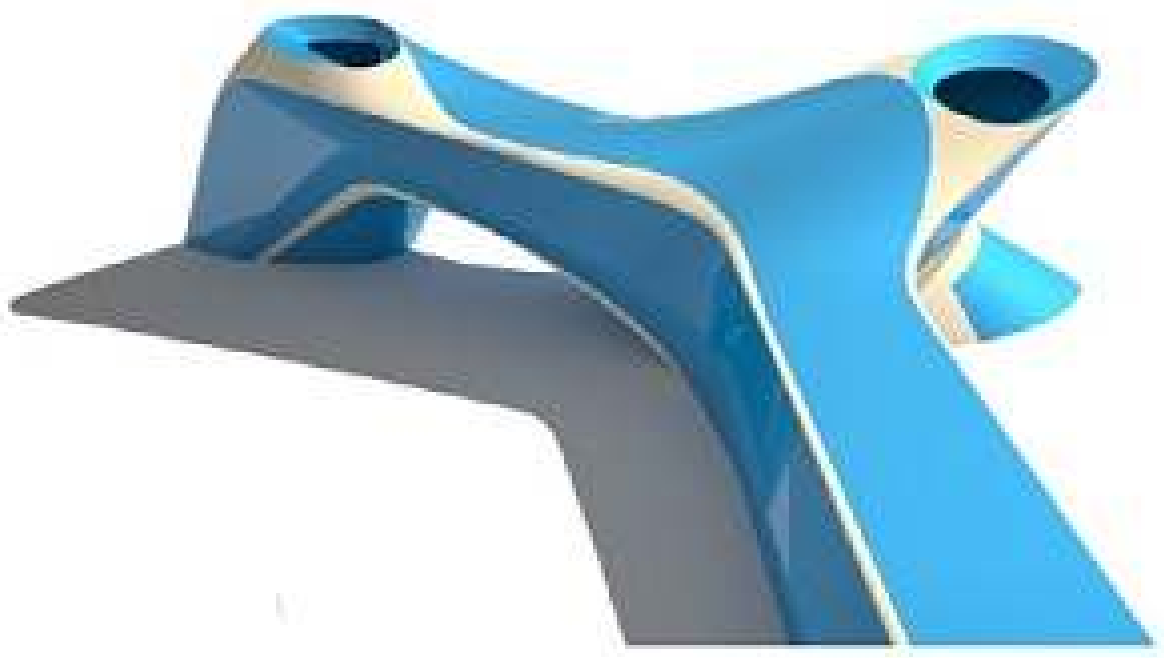

**Fig. 1.** Model of a shape represented by multiple strips of ruled surfaces [25].

In recent years, fundamental research has focused on the reverse problem or backward analysis: given a 3D curve, how can we characterize those surfaces that possess this curve as a special curve, rather than finding and classifying curves on analytical curved surfaces. Wang et al. [27] studied the problem of constructing a surface family from a given spatial geodesic. Kasap et al. [28] generalized the marching-scale functions of Wang and gave a sufficient condition for a given curve to be a geodesic on a surface. Lie et al. [29] derived the necessary and sufficient condition for a given curve to be the line of curvature on a surface. Bayram et al. [30] studied parametric surfaces which possess a given curve as a common asymptotic. However, they solved the problem using Frenet frame of the given curve.

In this paper, we obtain the necessary and sufficient condition for a given curve to be both isoparametric and asymptotic on a parametric surface depending on the RMF. Furthermore, we show that there exists ruled surfaces possessing a given curve as a common asymptotic curve and prove that there exists no developable surface passing through a given curve as an asymptotic curve except plane. We only study curves with an arc length parameter because such a study is easy to follow; if necessary, one can obtain similar results for arbitrarily parameterised regular curves.

## 2. Backgrounds

A parametric curve $r(s)$, $L_1 \le s \le L_2$, is a curve on a surface $P = P(s,t)$ in $\mathbb{R}^3$ that has a constant $s$ or $t$-parameter value. In this paper, $r'(s)$ denotes the derivative of $r$ with respect to arc length parameter $s$ and we assume that $r(s)$ is a regular curve, i.e. $r'(s) \neq 0$. For every point of $r(s)$, if $r''(s) \neq 0$, the set $\{T(s), N(s), B(s)\}$ is called the Frenet frame along $r(s)$, where $T(s) = r'(s)$, $N(s) = r''(s)/|r''(s)|$ and $B(s) = T(s) \times N(s)$ are the unit tangent, principal normal, and binormal vectors of the curve at the point $r(s)$, respectively. Derivative formulas of the Frenet frame is governed by the relations

$$\frac{d}{ds}\begin{pmatrix} T(s) \\ N(s) \\ B(s) \end{pmatrix} = \begin{pmatrix} 0 & \kappa(s) & 0 \\ -\kappa(s) & 0 & \tau(s) \\ 0 & -\tau(s) & 0 \end{pmatrix} \begin{pmatrix} T(s) \\ N(s) \\ B(s) \end{pmatrix} \tag{2.1}$$

where $\kappa(s) = \|r''(s)\|$ and $\tau(s) = \dfrac{(r'(s), r''(s), r'''(s))}{\|r'(s) \times r''(s)\|}$ are called the curvature and torsion of the curve $r(s)$, respectively [26].

Another useful frame along a curve is rotation minimizing frame (RMF). They are useful in animation, motion planning, swept surface constructions and related applications where the Frenet frame may prove unsuitable or undefined. A frame $\{T(s), U(s), V(s)\}$ among the frames on the curve $r(s)$ is called *rotation minimizing* if it is the frame of minimum twist around the tangent vector $T$. $\{T(s), U(s), V(s)\}$ is an RMF if

$$\begin{cases} U'(s) = -(U(s)\bullet r''(s))\, r'(s), \\ V'(s) = -(V(s)\bullet r''(s))\, r'(s), \end{cases} \tag{2.2}$$

where "•" denotes the standard inner product in $\mathbb{R}^3$ [32]. Observe that such a pair $U$ and $V$ is not unique; there exist a one parameter family of RMF's corresponding to different sets of initial values of $U$ and $V$. According to Bishop [20], a frame is an RMF if and only if each of $U'(s)$ and $V'(s)$ is parallel to $T(s)$. Equivalently,

$$\begin{cases} U'(s)\bullet V(s) \equiv 0, \\ V'(s)\bullet U(s) \equiv 0 \end{cases} \tag{2.3}$$

is the necessary and sufficient condition for the frame to be rotation minimizing [33]. Note that $\{T(s), U(s), V(s)\}$ is defined along the curve $r(s)$ even if the curvature vanishes where the Frenet frame is undefined.

There is a relation between the Frenet frame (if the Frenet frame is defined) and RMF, that is, $U$ and $V$ are the rotation of $N$ and $B$ of the curve $r(s)$ in the normal plane. Then,

$$\begin{bmatrix} U \\ V \end{bmatrix} = \begin{bmatrix} \cos\theta & \sin\theta \\ -\sin\theta & \cos\theta \end{bmatrix} \begin{bmatrix} N \\ B \end{bmatrix}, \tag{2.4}$$

where $\theta = \theta(s)$ is the angle between the vectors $N$ and $U$ (see Fig. 2), where $\theta' = -\tau$ [34].

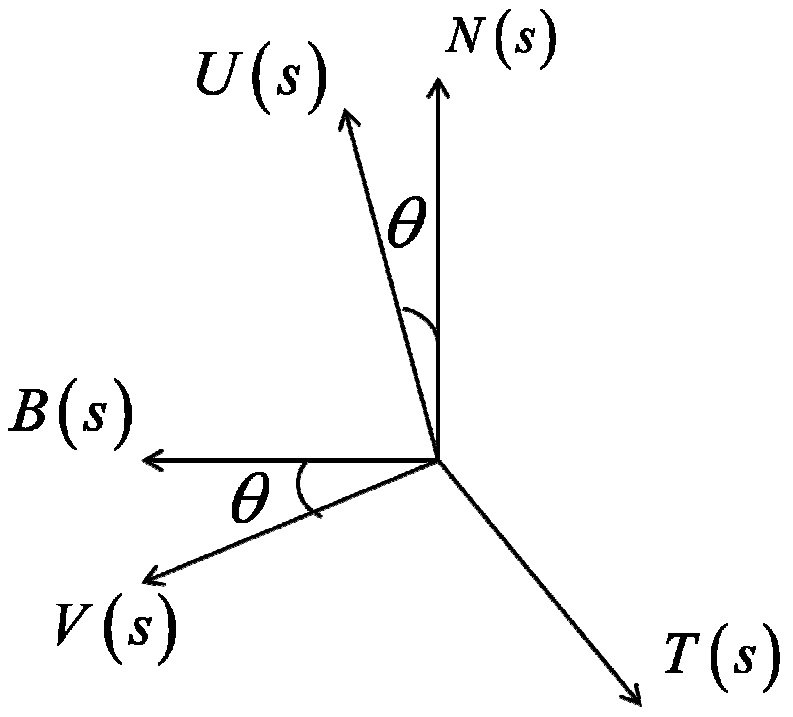

**Fig. 2** The Frenet frame *{T(s), N(s), B(s)}* and the vectors *U(s), V(s).*

Eqn. (2.4) implies that if $\{T(s), U(s), V(s)\}$ is an RMF then, it satisfies the following relations

$$U' = -\kappa\cos\theta T,\ V' = \kappa\sin\theta T,\ \theta' = -\tau. \tag{2.5}$$

## 3. Surface pencil with a given asymptotic curve

Suppose we are given a 3-dimensional parametric curve $r(s),\ L_1 \le s \le L_2,$ in which $s$ is the arc length (regular and $\|r'(\mathrm{s})\| = 1,\ L_1 \le s \le L_2$).

Surface pencil that interpolates $r(s)$ as a common curve is given in the parametric form as

$$P(s,t) = r(s) + a(s,t)T(s) + b(s,t)U(s) + c(s,t)V(s),\ L_1 \le s \le L_2,\ T_1 \le t \le T_2, \tag{3.1}$$

where $a(s,t), b(s,t)$ and $\mathrm{c}(s,t)$ are $C^1$ functions. The values of the *marching-scale functions* $a(s,t), b(s,t)$ and $\mathrm{c}(s,t)$ indicate, respectively, the extension-like, flexion-like, and retortion-like effects caused by the point unit through time *t,* starting from $r(s)$.

**Remark 3.1.** Observe that choosing different marching-scale functions yields different surfaces possessing $r(s)$ as a common curve.

Our goal is to find the necessary and sufficient conditions for which the curve $r(s)$ is isoparametric and asymptotic on the surface $P(s,t)$. Firstly, as $r(s)$ is an isoparametric curve on the surface $P(s,t)$, there exists a parameter $t_0 \in [T_1, T_2]$ such that

$$a(s,t_0) = b(s,t_0) = c(s,t_0) \equiv 0,\ L_1 \le s \le L_2,\ T_1 \le t_0 \le T_2. \tag{3.2}$$

Secondly the curve $r(s)$ is an asymptotic curve on the surface $P(s,t)$ if and only if

$$\frac{\partial n}{\partial s}(s,t_0) \bullet T(s) = 0, \tag{3.3}$$

where $n(s,t_0)$ is the surface normal along the curve $r(s)$. The normal vector of $P = P(s,t)$ can be written as

$$n(s,t) = \frac{\partial P(s,t)}{\partial s} \times \frac{\partial P(s,t)}{\partial t}.$$

From Eqns. (2.1) and (2.4), the normal vector can be expressed as

$$n(s,t)=\left(\frac{\partial c(s,t)}{\partial t}\left(a(s,t)\kappa(s)\cos\theta(s)+\frac{\partial b(s,t)}{\partial s}\right)-\frac{\partial b(s,t)}{\partial t}\left(\frac{\partial c(s,t)}{\partial s}-a(s,t)\kappa(s)\sin\theta(s)\right)\right)T(s)$$
$$+\left(\frac{\partial a(s,t)}{\partial t}\left(\frac{\partial c(s,t)}{\partial s}-a(s,t)\kappa(s)\sin\theta(s)\right)-\frac{\partial c(s,t)}{\partial t}\left(1+\frac{\partial a(s,t)}{\partial s}-b(s,t)\kappa(s)\cos\theta(s)+c(s,t)\kappa(s)\sin\theta(s)\right)\right)U(s)$$
$$+\left(\frac{\partial b(s,t)}{\partial t}\left(1+\frac{\partial a(s,t)}{\partial s}-b(s,t)\kappa(s)\cos\theta(s)+c(s,t)\kappa(s)\sin\theta(s)\right)-\frac{\partial a(s,t)}{\partial t}\left(a(s,t)\kappa(s)\cos\theta(s)+\frac{\partial b(s,t)}{\partial s}\right)\right)V(s).$$

Thus,

$$n(s,t_0)=\phi_1(s,t_0)T(s)+\phi_2(s,t_0)U(s)+\phi_3(s,t_0)V(s),$$

where

$$\phi_1(s,t_0)=\left(\frac{\partial c}{\partial t}(s,t_0)\left(a(s,t_0)\kappa(s)\cos\theta(s)+\frac{\partial b}{\partial s}(s,t_0)\right)-\frac{\partial b}{\partial t}(s,t_0)\left(\frac{\partial c}{\partial s}(s,t_0)-a(s,t_0)\kappa(s)\sin\theta(s)\right)\right),$$
$$\phi_2(s,t_0)=\left(\frac{\partial a}{\partial t}(s,t_0)\left(\frac{\partial c}{\partial s}(s,t_0)-a(s,t_0)\kappa(s)\sin\theta(s)\right)-\frac{\partial c}{\partial t}(s,t_0)\left(1+\frac{\partial a}{\partial s}(s,t_0)-b(s,t_0)\kappa(s)\cos\theta(s)+c(s,t_0)\kappa(s)\sin\theta(s)\right)\right),$$
$$\phi_3(s,t_0)=\left(\frac{\partial b}{\partial t}(s,t_0)\left(1+\frac{\partial a}{\partial s}(s,t_0)-b(s,t_0)\kappa(s)\cos\theta(s)+c(s,t_0)\kappa(s)\sin\theta(s)\right)-\frac{\partial a}{\partial t}(s,t_0)\left(a(s,t_0)\kappa(s)\cos\theta(s)+\frac{\partial b}{\partial s}(s,t_0)\right)\right).$$

**Remark 3.2.** Because of Eqn. (2.2) and the definition of partial differentiation we have

$$\frac{\partial a}{\partial s}(s,t_0)=\frac{\partial b}{\partial s}(s,t_0)=\frac{\partial c}{\partial s}(s,t_0)\equiv 0,\ \ t_0\in[T_1,T_2],\ L_1\le s\le L_2.$$

According to Remark 3.2, we should have

$$\phi_1(s,t_0)\equiv 0,\ \ \phi_2(s,t_0)=-\frac{\partial c}{\partial t}(s,t_0),\ \ \phi_3(s,t_0)=\frac{\partial b}{\partial t}(s,t_0).$$

From Eqns. (3.3) and (2.5), we get

$$\begin{aligned}\frac{\partial \hat{n}}{\partial s}(s,t_0)\cdot T(s)=0 &\Leftrightarrow \left(\frac{\partial\left(\phi_2(s,t)\mathrm{U}(s)+\phi_3(s,t)\mathrm{V}(s)\right)}{\partial s}\right)_{(s,t_0)}\cdot T(s)=0\\ &\Leftrightarrow \kappa\cos\theta\frac{\partial c}{\partial t}(s,t_0)=-\kappa\sin\theta\frac{\partial b}{\partial t}(s,t_0).\end{aligned}\tag{3.4}$$

Thus, we can present the following theorem:

**Theorem 3.3.** The necessary and sufficient condition for the curve $r(s)$ to be both isoparametric and asymptotic on the surface $P(s,t)$ is

$$\begin{cases} a(s,t_0)=b(s,t_0)=c(s,t_0)\equiv 0,\\ \kappa\cos\theta\dfrac{\partial c}{\partial t}(s,t_0)=-\kappa\sin\theta\dfrac{\partial b}{\partial t}(s,t_0),\\ \theta'(s)=-\tau(s).\end{cases}\tag{3.5}$$

**Corollary 3.4.** The sufficient condition for the curve $r(s)$ to be both isoparametric and asymptotic on the surface $P(s,t)$ is

$$
\begin{cases}
a(s,t_0) = f(s,t_0) \equiv 0, \\
b(s,t) = -f(s,t)\cos\theta,\ c(s,t) = f(s,t)\sin\theta, \\
\theta'(s) = -\tau(s).
\end{cases} \tag{3.6}
$$

## 4. Ruled surface pencil with a common asymptotic curve

**Theorem 4.1.** Given an arc length curve $r(s)$, there exists a ruled surface $P(s,t)$ possessing $r(s)$ as an asymptotic curve.

**Proof.** Choosing marching-scale functions as

$$
a(s,t) = (t-t_0)g(s),\ b(s,t) = (t_0 - t)\cos\theta,\ c(s,t) = (t-t_0)\sin\theta \tag{3.7}
$$

and $\theta'(s) = -\tau(s)$, Eqn. (3.1) takes the following form of a ruled surface

$$
P(s,t) = r(s) + (t-t_0)\left[g(s)T(s) - \cos\theta(s)U(s) + \sin\theta(s)V(s)\right], \tag{3.8}
$$

which satisfies Eqn. (3.6) interpolating $r(s)$ as an asymptotic curve.

**Corollary 4.2.** Ruled surface (3.8) is developable if and only if $r(s)$ is a planar curve.

**Proof.** $P(s,t) = r(s) + (t-t_0)\left[g(s)T(s) - \cos\theta(s)U(s) + \sin\theta(s)V(s)\right]$ is developable if and only if $\left(r', d, d'\right) = 0$, where $d(s) = g(s)T(s) - \cos\theta(s)U(s) + \sin\theta(s)V(s)$. Using Eqns. (2.1), (2.4) and (2.5) gives

$$
\begin{aligned}
d' &= g'T + gT' + \theta'\sin\theta U - \cos\theta U' + \theta'\cos\theta V + \sin\theta V' \\
&= g'T + g\kappa(\cos\theta U - \sin\theta V) - \tau\sin\theta U + \kappa\cos^2\theta T - \tau\cos\theta V + \kappa\sin^2\theta T \\
&= (g' + \kappa)T + (g\kappa\cos\theta - \tau\sin\theta)U - (g\kappa\sin\theta + \tau\cos\theta)V.
\end{aligned} \tag{3.9}
$$

Employing Eqn. (3.9) in the determinant we get $\tau = 0$, which implies that the curve $r(s)$ is planar.

**Corollary 4.3.** There exists no developable surface possessing a given curve as an asymptotic curve except plane.

**Proof.** Let $r(s)$ be any arc length curve. Because of Corollary 4.2 $r(s)$ should be planar. The surface will be in the form of Eqn. (3.8) and we have $a(s,t) = (t-t_0)g(s)$, $b(s,t) = (t_0 - t)\cos\theta$, $c(s,t) = (t-t_0)\sin\theta$. Thus, the surface normal is

$$
\begin{aligned}
n(s,t) = &\left(\frac{\partial a(s,t)}{\partial t}\left(\frac{\partial c(s,t)}{\partial s} - a(s,t)\kappa(s)\sin\theta(s)\right) - \frac{\partial c(s,t)}{\partial t}\left(1 + \frac{\partial a(s,t)}{\partial s} - b(s,t)\kappa(s)\cos\theta(s) + c(s,t)\kappa(s)\sin\theta(s)\right)\right)U(s) \\
&+ \left(\frac{\partial b(s,t)}{\partial t}\left(1 + \frac{\partial a(s,t)}{\partial s} - b(s,t)\kappa(s)\cos\theta(s) + c(s,t)\kappa(s)\sin\theta(s)\right) - \frac{\partial a(s,t)}{\partial t}\left(a(s,t)\kappa(s)\cos\theta(s) + \frac{\partial b(s,t)}{\partial s}\right)\right)V(s) \\
&= \left((t_0 - t)g^2\kappa\sin\theta - \sin\theta\left(1 + (t-t_0)\kappa\cos^2\theta + (t-t_0)\kappa\sin^2\theta\right)\right)U + \left(-\cos\theta\left(1 + (t-t_0)\kappa\cos^2\theta + (t-t_0)\kappa\sin^2\theta\right) - g\left((t-t_0)g\kappa\cos\theta\right)\right)V \\
&= -\sin\theta\left((t-t_0)g^2\kappa + 1 + (t-t_0)\kappa\right)U - \cos\theta\left((t-t_0)g^2\kappa + 1 + (t-t_0)\kappa\right)V.
\end{aligned}
$$

Using Eqn. (2.4) the surface normal takes the form $n(s,t) = -\left((t-t_0)g^2\kappa + 1 + (t-t_0)\kappa\right)B(s)$ which is always parallel to the binormal $B(s)$ of the curve $r(s)$. Since $r(s)$ is planar, $B(s)$ is constant giving that the surface $P(s,t)$ is a part of a plane.

## 5. Examples of generating surfaces with a common asymptotic curve

**Example 5.1.** Let $r(s)=(\cos(s),\sin(s),0)$ be a unit speed curve. It is obvious that

$$T(s)=\left(-\sin(s),\cos(s),0\right),\ \kappa(s)=1,\ \tau(s)=0.$$

If we take $U(s)=(0,0,1)$ and $V(s)=(\cos s,\sin s,0)$, then Eqn. (2.2) is satisfied and $\{T(s),U(s),V(s)\}$ is an RMF. By choosing marching-scale functions as $a(s,t)\equiv 0,\ b(s,t)=\sin t-1,\ c(s,t)=\cos t$ and $t_0=\theta=\frac{\pi}{2}$, then Eqn. (3.5) is satisfied. Thus, we immediately obtain a member of the surface pencil with a common asymptotic curve $r(s)$ as

$$P_1(s,t)=\left(\cos(s)(\cos t+1),\ \sin(s)(\cos t+1),\ \sin t-1\right),$$

where $0\le s\le 2\pi,\ 0\le t\le 2\pi$ (Fig. 2).

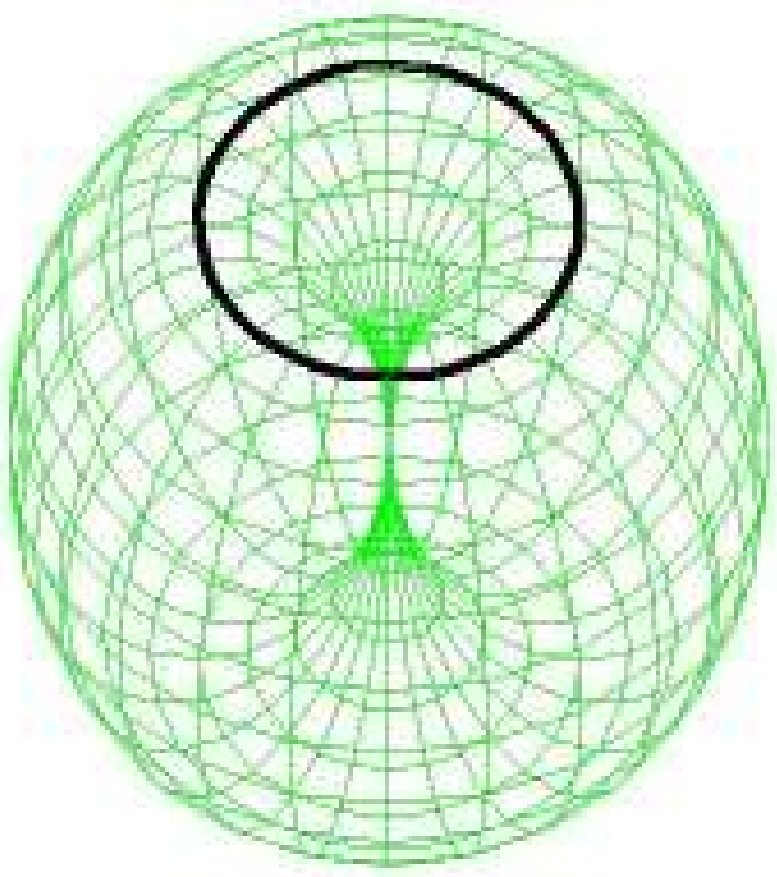

**Fig. 2.** $P_1(s,t)$ as a member of the surface pencil and its asymptotic curve.

**Example 5.2.** Let $r(s)=\left(\frac{3}{5}\sin(s),\frac{3}{5}\cos(s),\frac{4}{5}s\right)$ be a unit speed curve. Then it is easy to show that

$$\begin{cases} T(s)=\left(\frac{3}{5}\cos(s),-\frac{3}{5}\sin(s),\frac{4}{5}\right), \\ \kappa(s)=3/5,\ \tau(s)=-4/5. \end{cases}$$

If we choose $\theta=4s/5$ and

$$\begin{cases} U'=\left(-\frac{9}{25}\cos(s)\cos\left(\frac{4s}{5}\right),\frac{9}{25}\sin(s)\cos\left(\frac{4s}{5}\right),-\frac{12}{25}\cos\left(\frac{4s}{5}\right)\right), \\ V'=\left(\frac{9}{25}\cos(s)\sin\left(\frac{4s}{5}\right),-\frac{9}{25}\sin(s)\sin\left(\frac{4s}{5}\right),\frac{12}{25}\sin\left(\frac{4s}{5}\right)\right), \end{cases}$$

then Eqn. (2.5) is satisfied. By integration, we obtain

$$\begin{cases} U=\left(-\frac{1}{10}\sin\left(\frac{9s}{5}\right)-\frac{9}{10}\sin\left(\frac{s}{5}\right),\ -\frac{1}{10}\cos\left(\frac{9s}{5}\right)-\frac{9}{10}\cos\left(\frac{s}{5}\right),\ -\frac{3}{5}\sin\left(\frac{4s}{5}\right)\right), \\ V=\left(\frac{9}{10}\cos\left(\frac{s}{5}\right)-\frac{1}{10}\cos\left(\frac{9s}{5}\right),\ \frac{1}{10}\sin\left(\frac{9s}{5}\right)-\frac{9}{10}\sin\left(\frac{s}{5}\right),\ -\frac{3}{5}\cos\left(\frac{4s}{5}\right)\right). \end{cases}$$

Now, $\{T(s), U(s), V(s)\}$ is an RMF since it satisfies Eqn. (2.2). If we take $a(s,t)\equiv 0$, $b(s,t)=e^{s}t\cos\left(\frac{4s}{5}\right)$, $c(s,t)=-e^{s}t\sin\left(\frac{4s}{5}\right)$ and $t_0=0$, then Eqn. (3.6) is satisfied. Thus, we obtain a member of the surface pencil with a common asymptotic curve $r(s)$ as

$$P_2(s,t)=\left(\frac{3}{5}\sin(s)+e^{s}t\left(\cos\left(\frac{4s}{5}\right)\left(-\frac{1}{10}\sin\left(\frac{9s}{5}\right)-\frac{9}{10}\sin\left(\frac{s}{5}\right)\right)-\sin\left(\frac{4s}{5}\right)\left(-\frac{1}{10}\cos\left(\frac{9s}{5}\right)+\frac{9}{10}\cos\left(\frac{s}{5}\right)\right)\right),\right.$$
$$\left.\frac{3}{5}\cos(s)+e^{s}t\left(\cos\left(\frac{4s}{5}\right)\left(-\frac{1}{10}\cos\left(\frac{9s}{5}\right)-\frac{9}{10}\cos\left(\frac{4s}{5}\right)\right)-\sin\left(\frac{4s}{5}\right)\left(-\frac{9}{10}\sin\left(\frac{s}{5}\right)+\frac{1}{10}\sin\left(\frac{9s}{5}\right)\right)\right),\frac{4}{5}s\right),$$

where $-\pi\le s\le 0$, $-1\le t\le 1/2$ (Fig. 3).

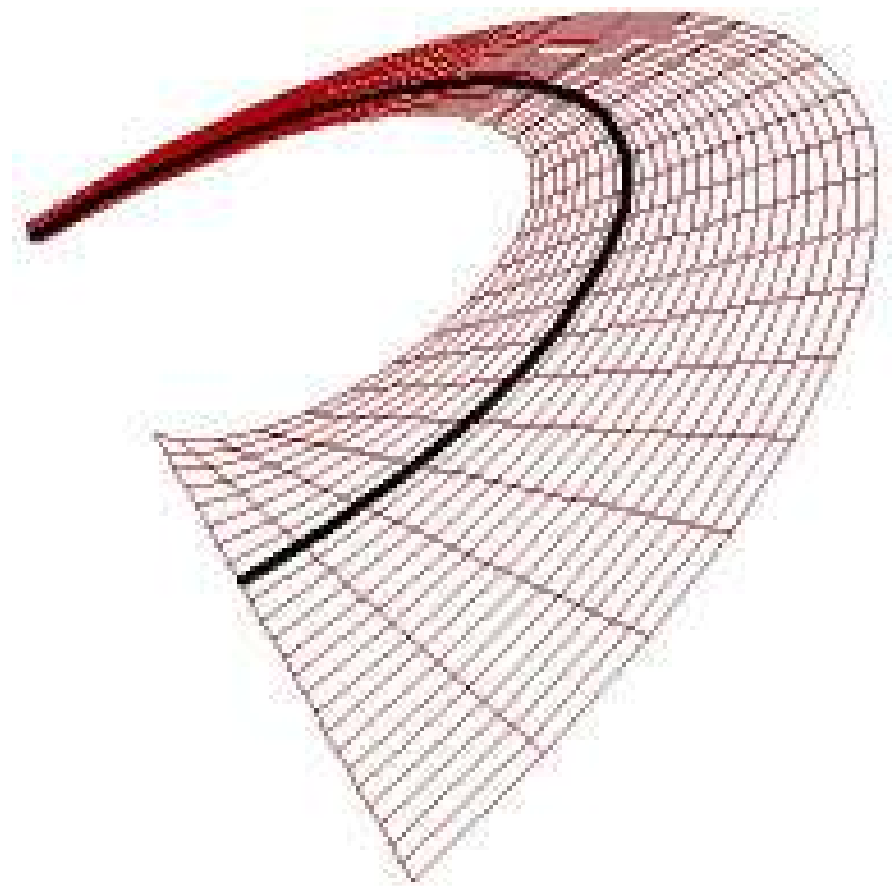

**Fig. 3.** $P_2(s,t)$ as a member of the surface pencil and its asymptotic curve.

For the same curve, if we take marching-scale functions as $a(s,t)\equiv 0$, $b(s,t)=-\cos(\frac{4s}{5})$, $c(s,t)=\sin\left(\frac{4s}{5}\right)$ and $t_0=0$ then by Eqn. (3.8) we obtain the following ruled surface with a common asymptotic curve $r(s)$ as

$$P_3(s,t)=\left(\frac{3}{5}\sin(s)+\frac{t}{10}\left[\cos\left(\frac{4s}{5}\right)\left[\sin\left(\frac{9s}{5}\right)+9\sin\left(\frac{s}{5}\right)\right]+\sin\left(\frac{4s}{5}\right)\left[9\cos\left(\frac{s}{5}\right)-\cos\left(\frac{9s}{5}\right)\right]\right],\right.$$
$$\left.\frac{3}{5}\cos(s)+\frac{t}{10}\left[\cos\left(\frac{4s}{5}\right)\left[\cos\left(\frac{9s}{5}\right)+9\cos\left(\frac{s}{5}\right)\right]+\sin\left(\frac{4s}{5}\right)\left[\sin\left(\frac{9s}{5}\right)-9\sin\left(\frac{s}{5}\right)\right]\right],\frac{4}{5}s\right),$$

where $-\pi\le s\le 0$, $-1\le t\le \frac{1}{2}$ (Fig. 4).

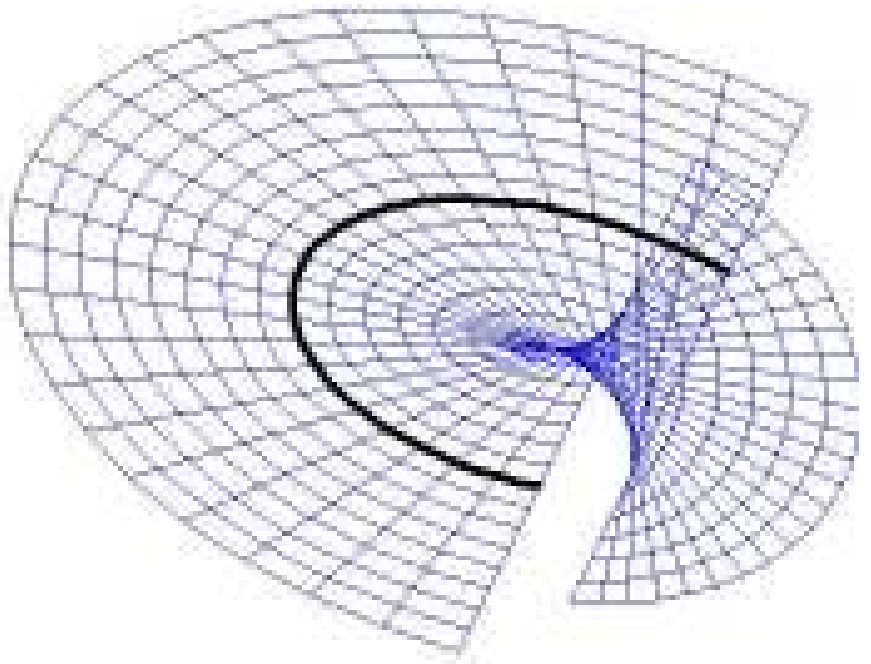

**Fig. 4.** Ruled surface $P_3(s,t)$ as a member of the surface pencil and its asymptotic curve.

**Example 5.3.** Let $r(s)=\left(\frac{\sqrt{3}}{2}\sin(s),\frac{s}{2},\frac{\sqrt{3}}{2}\cos(s)\right)$ be a unit speed curve. Then,

$$T(s)=\left(\frac{\sqrt{3}}{2}\cos(s),\frac{1}{2},-\frac{\sqrt{3}}{2}\sin(s)\right),\ \kappa(s)=\frac{\sqrt{3}}{2},\ \tau(s)=-\frac{1}{2}.$$

Using Eqns. (2.2) and (2.5) we obtain $\theta=-s/2$ and

$$\begin{cases} U(s)=\left(-\frac{1}{4}\sin\left(\frac{3s}{2}\right)-\frac{3}{4}\sin\left(\frac{s}{2}\right),-\frac{\sqrt{3}}{2}\sin\left(\frac{s}{2}\right),-\frac{1}{4}\cos\left(\frac{3s}{2}\right)-\frac{3}{4}\cos\left(\frac{s}{2}\right)\right), \\ V(s)=\left(\frac{1}{4}\cos\left(\frac{3s}{2}\right)-\frac{3}{4}\cos\left(\frac{s}{2}\right),\frac{\sqrt{3}}{2}\cos\left(\frac{s}{2}\right),\frac{3}{4}\sin\left(\frac{s}{2}\right)-\frac{1}{4}\sin\left(\frac{3s}{2}\right)\right). \end{cases}$$

Now, $\{T(s), U(s), V(s)\}$ is an RMF. If we take $a(s,t)\equiv 0,\ b(s,t)=s^2t,$ $c(s,t)=s^2t\tan\left(\frac{s}{2}\right)$ and $t_0=0$, then Eqn. (3.5) is satisfied. Thus, we have a member of the surface pencil with a common asymptotic curve $r(s)$ as

$$\begin{aligned} P_4(s,t)=\Bigg(&\frac{\sqrt{3}}{2}\sin(s)-\frac{s^2t}{4}\left[\sin\left(\frac{3s}{2}\right)+3\sin\left(\frac{s}{2}\right)-\tan\left(\frac{s}{2}\right)\left(\cos\left(\frac{3s}{2}\right)+3\cos\left(\frac{s}{2}\right)\right)\right], \\ &\frac{s}{2}-s^2t\left[\frac{\sqrt{3}}{2}\sin\left(\frac{s}{2}\right)-\frac{\sqrt{3}}{2}\tan\left(\frac{s}{2}\right)\cos\left(\frac{s}{2}\right)\right], \\ &\frac{\sqrt{3}}{2}\cos(s)-\frac{s^2t}{4}\left[\cos\left(\frac{3s}{2}\right)+3\cos\left(\frac{s}{2}\right)-\tan\left(\frac{s}{2}\right)\left(-\sin\left(\frac{3s}{2}\right)+3\sin\left(\frac{s}{2}\right)\right)\right]\Bigg), \end{aligned}$$

where $-2<s\leq 2$, $-1\leq t\leq 1$ (Fig. 5).

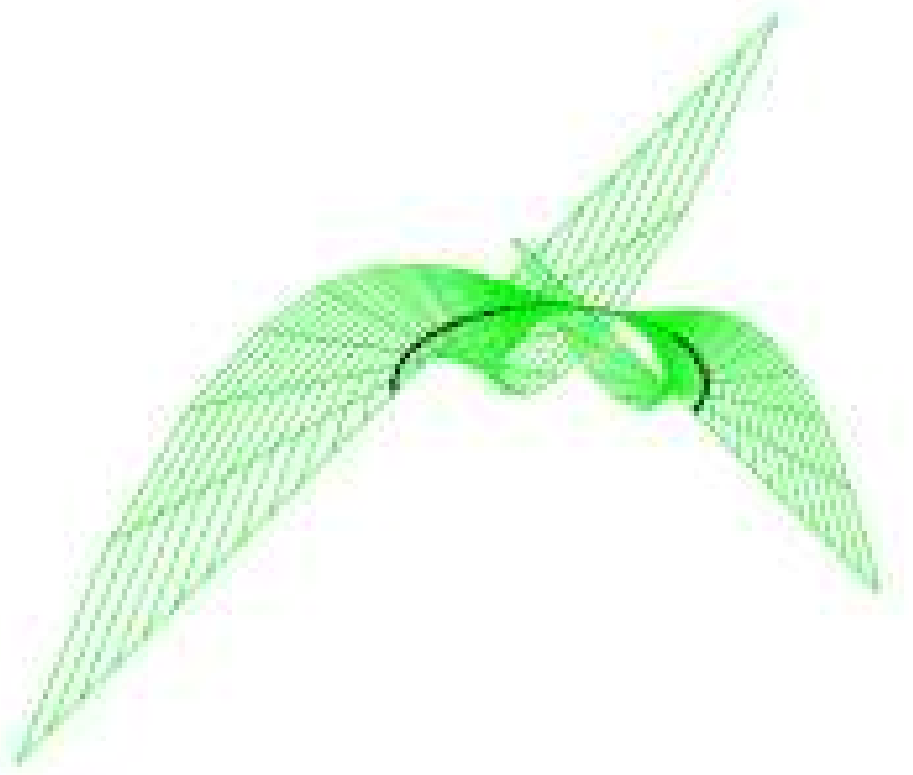

**Fig. 5.** $P_4(s,t)$ as a member of the surface pencil and its asymptotic curve.

For the same curve let us find a ruled surface. In Eqn. (3.8), if we take $g(s)=\cos(s)$, then we obtain the following ruled surface with a common asymptotic curve $r(s)$ as

$$P_5(s,t)=\Bigg(\frac{\sqrt{3}}{2}\sin(s)+t\left(\frac{\sqrt{3}}{2}\cos^2(s)+\frac{1}{4}\sin\left(\frac{3s}{2}\right)\cos\left(\frac{s}{2}\right)+\frac{3}{2}\sin\left(\frac{s}{2}\right)\cos\left(\frac{s}{2}\right)-\frac{1}{4}\cos\left(\frac{3s}{2}\right)\sin\left(\frac{s}{2}\right)\right),$$
$$\frac{s}{2}+t\left(\frac{1}{2}\cos(s)\right),$$
$$\frac{\sqrt{3}}{2}\cos(s)+t\left(-\frac{\sqrt{3}}{2}\sin(s)\cos(s)+\frac{1}{4}\cos\left(\frac{3s}{2}\right)\cos\left(\frac{s}{2}\right)+\frac{1}{4}\sin\left(\frac{3s}{2}\right)\sin\left(\frac{s}{2}\right)+\frac{3}{4}\cos(s)\right)\Bigg),$$

where $-2<s\leq 2$, $-1\leq t\leq 1$ (Fig. 6).

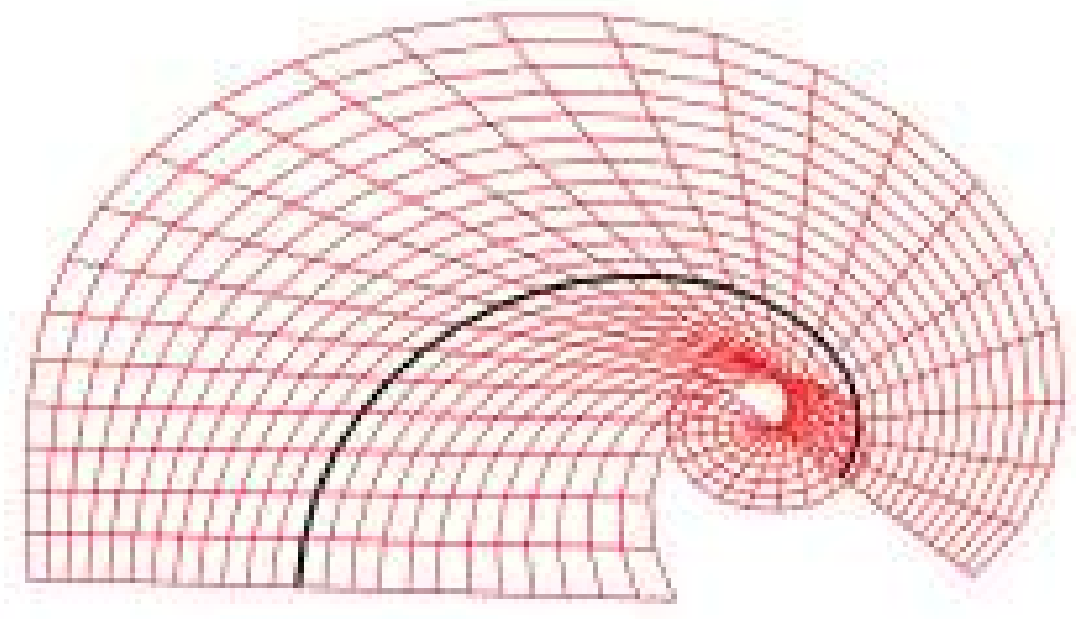

**Fig. 6.** $P_5(s,t)$ as a member of the ruled surface pencil and its asymptotic curve.

If we take $g(s)=\cos(s)\sin(s)$ then, we obtain the ruled surface with a common asymptotic curve $r(s)$ as

$$P_6(s,t)=\left(\frac{\sqrt{3}}{2}\sin(s)+t\left(\frac{\sqrt{3}}{2}\cos(s)\sin(s)+\frac{1}{4}\sin\left(\frac{3s}{2}\right)\cos\left(\frac{s}{2}\right)+\frac{3}{2}\sin\left(\frac{s}{2}\right)\cos\left(\frac{s}{2}\right)-\frac{1}{4}\cos\left(\frac{3s}{2}\right)\sin\left(\frac{s}{2}\right)\right),\right.$$

$$\frac{s}{2}+t\left(\frac{1}{2}\cos(s)\sin(s)\right),$$

$$\left.\frac{\sqrt{3}}{2}\cos(s)+t\left(-\frac{\sqrt{3}}{2}\sin^2(s)\cos(s)+\frac{1}{4}\cos\left(\frac{3s}{2}\right)\cos\left(\frac{s}{2}\right)+\frac{1}{4}\sin\left(\frac{3s}{2}\right)\sin\left(\frac{s}{2}\right)+\frac{3}{4}\cos(s)\right)\right),$$

where $-2 < s \leq 2$, $-1 \leq t \leq 1$ (Fig 7).

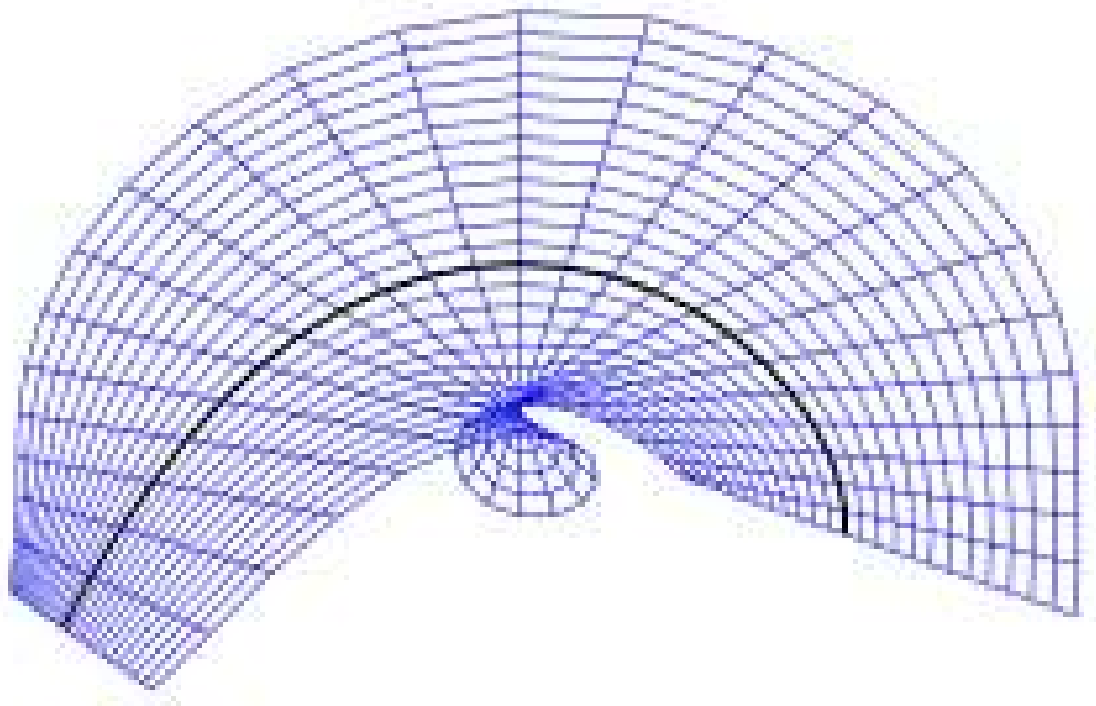

**Fig. 7.** $P_6(s,t)$ as a member of the ruled surface pencil and its asymptotic curve.

## 6. Conclusion

In this paper, we introduce a method for finding a surface pencil passing through a given asymptotic curve as an isoparametric curve. It has been shown that if the necessary and sufficient conditions are satisfied (i.e., Eqn. 3.5) for a given curve, a common isoparametric curve and asymptotic curve on a surface pencil can always be found.

There are several opportunities for further work. One possible alternative is to consider the realm of implicit surfaces $F(x,y,z)=0$ and to attempt to find the constraints for a given curve $r(s)$ to be asymptotic on $F(x,y,z)=0$. In addition, an anologue of the problem addressed in this paper may be considered for 3-surfaces in 4-space or other types of marching-scale functions may be investigated.

## 7. Acknowledgements

Ergin Bayram would like to thank TUBITAK (The Scientific and Technological Research Council of Turkey) for their financial supports during his doctorate studies.